\documentclass[11 pt] {amsart}

\usepackage[margin=1in]{geometry}
\usepackage{graphicx}
\usepackage{latexsym}
\usepackage{amsfonts}
\usepackage{amsmath}
\usepackage{enumerate}
\usepackage{hyperref}
\usepackage{comment}

\newcommand\blackslug{\hbox{\hskip 1pt \vrule width 4pt height 8pt depth 1.5pt
        \hskip 1pt}}
\newcommand\bbox{\hfill \quad \blackslug \medbreak}
\def\d{\hbox{-}}
\def\c{\hbox{-}\cdots\hbox{-}}

\newtheorem{theorem}{Theorem}[section]

\newtheorem{prop}[theorem]{Proposition}

\newtheorem{conj}[theorem]{Conjecture} 

\newcounter{claim}

\newcommand{\Proof}{\setcounter{claim}{0}\noindent{\bf Proof.}\ \ }

\title{Even pairs in Berge graphs with no balanced skew-partitions}

\author{Tara Abrishami$^{\ast \dagger}$}
\author{Maria Chudnovsky$^{\ast \mathparagraph}$}
\author{Yaqian Tang$^{\ast \amalg}$}
\address{$^{\ast}$Princeton University, Princeton, NJ, USA}
\address{$^{\dagger}$ Supported by NSF-EPSRC Grant DMS-2120644.}
\address{$^{\mathparagraph}$ Supported by NSF-EPSRC Grant DMS-2120644 and by AFOSR grant FA9550-22-1-0083.}
\address{$^{\amalg}$ Supported by AFOSR grant FA9550-22-1-0083.}

\usepackage[foot]{amsaddr}
\begin{document}
\maketitle

\begin{abstract}
Let $G$ be a Berge graph that has no odd prism and no antihole of length at least six as an induced subgraph. We show that every such graph $G$ with no balanced skew-partition is either complete or has an even pair. 
\end{abstract}

\section{Introduction}
All graphs in this paper are finite and simple. Let $\chi(G)$ and $\omega(G)$ denote the chromatic number and the clique number of a graph $G$, respectively. A graph $G$ is \emph{perfect} if every induced subgraph $H$ of $G$ satisfies $\chi(H)=\omega(H)$. The {\em complement} of a graph $G$, denoted by $\overline{G}$, has the same vertex set as $G$, and two distinct vertices in $\overline{G}$ are adjacent if and only if they are not adjacent in $G$. A \emph{hole} in a graph $G$ is an induced subgraph isomorphic to a cycle on at least five vertices, and an \emph{antihole} is an induced subgraph whose complement is a hole in $\overline{G}$ (note that in this definition induced cycles of length four are treated in a non-standard way). The {\em length} of a hole (antihole) is equal to the number of its vertices. A graph is \emph{Berge} if it contains no odd hole and no odd antihole as an induced subgraph. In the 1960s, Berge \cite{Berge} conjectured that a graph is perfect if and only if it is \emph{Berge}. The study of perfect graphs became a major area of research in structural graph theory after Berge's conjecture. In 2002, Chudnovsky, Robertson, Seymour, and Thomas \cite{SPGT} proved the conjecture, which then became known as the \emph{Strong Perfect Graph Theorem (SPGT)}. 

An {\em even pair} in a graph $G$ is a pair $\{u,v\}$ of nonadjacent vertices such that every induced path from $u$ to $v$ in $G$ has an even number of edges. Before the SPGT was proved, many results focused on properties of {\em minimal imperfect graphs}: imperfect graphs $G$ such that every proper induced subgraph of $G$ is perfect. In particular, Meyniel \cite{meyniel} proved that minimal imperfect graphs do not have an even pair. Also, the proof of the SPGT was simplified by Chudnovsky and Seymour in 2007 using even pairs \cite{EPIBG}. 

A graph $G$ is {\em complete} if every pair of vertices in $G$ is adjacent. For a vertex $v\in V(G)$, we denote the set of vertices adjacent to $v$ by
$N_G(v)$ or by $N(v)$. We say a graph $G'$ is obtained by {\em contracting an even pair} $\{u,v\}$ in $G$ if:
\begin{itemize}
    \item $V(G')=(V(G)\setminus \{u,v\})\cup \{w\}$ (where $w \not \in V(G)$);
    \item $G'\setminus \{w\}=G\setminus \{u,v\}$; and
    \item $N_{G'}(w)=N_{G}(u)\cup N_{G}(v)$
\end{itemize}
We denote the graph obtained by contracting the even pair $\{u,v\}$ by $G/\{u,v\}$. A \emph{sequence of contraction} for a graph $G$ is a sequence of graphs $G_0,\cdots,G_k$ such that $G_0=G$, $G_k$ has no even pair, and for all $0 \leq i \leq k-1$, there exists an even pair $\{u, v\}$ in $G_i$ such that $G_{i+1}=G_{i}/\{u,v\}$. A graph is \emph{even-contractile} if it has a sequence of contraction with $G_k$ being a complete graph. Fonlupt and Uhry \cite{FONLUPT} observed that if $G$ is Berge with an even pair $\{u,v\}$, then $G/\{u,v\}$ is also Berge and $\omega(G/\{u,v\})=\omega(G)$. It follows that  given a $\chi(G/\{u,v\})$-coloring of $G/\{u,v\}$, one can obtain a $\chi(G)$-coloring of $G$ by preserving the same colors for vertices in $G\setminus \{u,v\}$ and assigning $u,v$ the color of the additional vertex $w$. Also, an even pair can be recognized in polynomial time using the algorithm for detecting an odd hole (\cite{recogoddhole} and \cite{threeinatree}): Given a pair of nonadjacent of vertices $\{u,v\}$, we add a new vertex $p$ to $G$ such that $N(p)=\{u,v\}$; the new graph
contains no odd hole if and only if $\{u,v\}$ is an even pair in $G$.
Now suppose that  we can find a sequence of contractions $G_0,\cdots, G_k$ for a Berge graph $G$ (in polynomial time), and such that 
the graph $G_k$ belongs to a simpler class of graphs, where a polynomial-time
coloring algorithm is known. In this case a $\chi(G)$-coloring of $G$ can be derived from a $\chi(G_k)$-coloring of $G_k$  in polynomial time. An example of such a situation is  when $G$ is even-contractile, as finding a $\chi(G_k)$-coloring of the complete graph $G_k$ is trivial. Therefore, a natural question is to identify which Berge graph is even-contractile.

To this end, Everett and Reed \cite{everett2001even} proposed a conjecture for characterizing even-contractile Berge graphs. A \emph{prism} in a graph $G$ is an induced subgraph that consists of two disjoint cliques $\{a_1,a_2,a_3 \}$ and $\{b_1,b_2,b_3\}$ and three disjoint paths $P_1,P_2,P_3$ from $a_i$ to $b_i$ for each $i$, and with no other edge except for those in the two cliques and in the three paths. The paths $P_1, P_2, P_3$ are called the {\em rungs} of the prism. Under these conditions, a prism is \emph{odd} ({\em even}) if all the three rungs have odd (even) number of edges. Note that every prism in a Berge graph is either even or odd. Everett and Reed conjectured the following:

\begin{conj}[\cite{everett2001even}]\label{Conj}
A Berge graph with no induced subgraph isomorphic to an antihole of length at least six or an odd prism is even-contractile.
\end{conj}

This conjecture remains open, but several related theorems have been proved. Maffray and Trotignon \cite{MAFFRAY} showed that a Berge graph that has no prism and no antihole of length at least six is even-contractile. Chudnovsky, Maffray, Seymour, and Spirkl \cite{EPSQ} showed that if a Berge graph contains has no cycle on four vertices and no odd prism of a particular type, then it is either complete or has an even pair. 

The main theorem of this paper is the following:

\begin{theorem}\label{thm:main}
Let $G$ be a Berge graph with no induced subgraph isomorphic to an antihole of length at least six or an odd prism. If $G$ does not admit a balanced skew-partition, then $G$ is either complete or has an even pair.
\end{theorem}

Note that Theorem~\ref{thm:main} does not imply Conjecture~\ref{Conj}
(even in the restricted graph class) since it only guarantees the existence
of one even pair, contracting which may introduce a skew-partition.
A {\em balanced skew-partition} is a type of decomposition that appears in the proof of the SPGT. In 2003, Chudnovsky \cite{thesis} proved a structural decomposition theorem for {\em trigraphs}, which is a generalization of graphs with possible ``undecided'' edges called {\em switchable pairs}. In particular, the theorem implies that a Berge graph either belongs to some ``basic'' class, or has a balanced skew-partition, or a $2$-join, or a $2$-join in the complement. Our result is based on this decomposition theorem, and the notion of trigraph is very helpful to the proof. 

The remainder of the paper is organized as follows. In Section \ref{trigraphs}, we introduce the definitions related to trigraphs and present relevant theorems that have been proved. We also define \emph{basic} trigraphs and decompositions, namely balanced skew-partition, $2$-join, and the complement of $2$-join. In Section \ref{decomp}, we define a class $\mathcal{F}$ of Berge trigraphs and a subclass called \emph{favorable trigraphs} that interact well with the $2$-join decomposition. In particular, we will show that almost all trigraphs in $\mathcal{F}$ are favorable when forbidding antihole of length six and balanced skew-partition. In Section \ref{Basic}, we show that basic trigraphs have even pairs, and favorable basic trigraphs have even pairs in certain desirable location. In Section \ref{BofD}, we apply the technique of \emph{block of decompositions} introduced in \cite{algoTri} to handle $2$-join and its complement. This technique allows us to decompose any trigraph in $\mathcal{F}$ with no balanced skew-partition into basic trigraphs while keeping track of even pairs. Finally, we prove a generalization of our main theorem \ref{thm:main} for trigraphs.

\section{Trigraphs}\label{trigraphs}
In this paper, we mainly adopt the notation regarding trigraphs from the work by Chudnovsky, Trotignon, Trunck, and Vušković \cite{algoTri}. For the sake of clarity, we restate relevant definitions and introduce new definitions that will appear in the paper. 

For a set $X$, we denote by $\binom{X}{2}$ the set of all subsets of $X$ of size $2$. For brevity, an element $\{u,v\}$ of $\binom{X}{2}$ is also denoted by $uv$, or equivalently, $vu$. A \emph{trigraph} $T$ consists of a finite  set $V(T)$, called the \emph{vertex set} of $T$, and a map $\theta: \binom{V(T)}{2}\to \{-1,0,1 \}$, called the \emph{adjacency function} of $T$. Two distinct vertices $u,v$ of $T$ are \emph{strongly adjacent} if $\theta(uv)=1$, \emph{strongly antiadjacent} if $\theta(uv)=-1$, and \emph{semiadjacent} if $\theta(uv)=0$. We say that $u$ and $v$ are \emph{adjacent} if $\theta(uv) \in \{0,1\}$ and \emph{antiadjacent} if $\theta(uv) \in \{0,-1\}$. If $u$ and $v$ are adjacent (antiadjacent), we also say that $u$ is \emph{adjacent} (\emph{antiadjacent}) to $v$, or $u$ is a \emph{neighbor} (\emph{antineighbor}) of $v$. Similarly, if $u$ and $v$ are strongly adjacent (strongly antiadjacent), we say $u$ is a \emph{strong neighbor} (\emph{strong antineighbor}) of $v$. For $v \in V(T)$, let $N(v)$ denote the set of all vertices in $V(T) \setminus \{v\}$ that are adjacent to $v$, and let $N[v]$ denote $N(v)\cup \{v\}$. An \emph{edge} (\emph{antiedge}) is a pair of adjacent (antiadjacent) vertices. A \emph{switchable pair} is a pair of semiadjacent vertices, and a \emph{strong edge} (\emph{antiedge}) is a pair of strongly adjacent (strongly antiadjacent) vertices. An edge $uv$ (antiedge, strong edge, strong antiedge, switchable pair) is \emph{between} two sets $A\subseteq V(T)$ and $B \subseteq V(T)$ if $u\in A$ and $v\in B$, or if $u\in B$ and $v\in A$.

Let $T$ be a trigraph. The \emph{complement} of $T$, denoted by $\overline{T}$, is a trigraph with $V(\overline{T})=V(T)$ and the adjacency function $\overline{\theta}=-\theta$. Let $A\subset V(T)$ and $b\in V(T) \setminus A$. We say that $b$ is \emph{strongly complete} (\emph{strongly anticomplete}) to $A$ if $b$ is strongly adjacent (strongly antiadjacent) to every vertex of $A$; $b$ is \emph{complete} (\emph{anticomplete}) to $A$ if $b$ is adjacent (antiadjacent) to every vertex of $A$. For two disjoint subsets $A \subset V(T)$ and $B\subset V(T)$, $B$ is \emph{strongly complete} (\emph{strongly anticomplete}, \emph{complete}, \emph{anticomplete}) to $A$ if every vertex of $B$ is strongly complete (strongly anticomplete, complete, anticomplete) to $A$. 

A \emph{clique} of $T$ is set of pairwise adjacent vertices of $T$, and a \emph{strong clique} is a set of pairwise strongly adjacent vertices of $T$. A trigraph $T$ is \emph{complete} if $V(T)$ is a clique. A \emph{stable set} of $T$ is a set of pairwise antiadjacent vertices of $T$. For $X \subseteq V(T)$, the trigraph \emph{induced by $T$ on $X$}, denoted by $T|X$, has vertex set $X$ and adjacency function $\theta|_{X}$, the restriction of $\theta$ to $\binom{X}{2}$. We denote by $T\setminus X$ the trigraph $T|(V(T)\setminus X)$. Isomorphism between trigraphs is defined in the natural way. For two trigraphs $T$ and $H$, $H$ is an \emph{induced subtrigraph} of $T$ (or \emph{T contains H as an induced subgtrigraph}) if $H$ is isomorphic to $T|X$ for some $X \subseteq V(T)$. Since this paper mainly considers the induced subtrigraph containment relation, we say that $T$ \emph{contains} $H$ if $T$ contains $H$ as an induced subtrigraph.

Let $\eta(T)$ denote the set of all strong edges of $T$, $\nu(T)$ the set of all strong antiedges of $T$, $\sigma(T)$ the set of all switchable pairs of $T$. If $\sigma(T)$ is empty, $T$ is a \emph{graph}. A \emph{semirealization} of $T$ is a trigraph $T'$ with vertex set $V(T)$ that satisfies $\eta(T) \subseteq \eta(T')$ and $\nu(T) \subseteq \nu(T')$. A \emph{realization} of $T$ is any graph that is semirealization of $T$. For $S\subseteq \sigma(T)$, we denote by $G_S^T$ the realization of $T$ with edge set $\eta(T) \cup S$. The realization $G^T_{\sigma(T)}$ is called the \emph{full realization} of $T$.

Let $T$ be a trigraph. For $X\subseteq V(T)$, we say that $X$ and $T|X$ are \emph{connected} (\emph{anticonnected}) if the graph $G_{\sigma(T|X)}^{T|X}$ ($\overline{G_{\emptyset}^{{T|X}}}$) is connected. A \emph{connected component} (or simply \emph{component}) of $X$ is maximal connected subset of $X$, and an \emph{anticonnected component} (or simply \emph{anticomponent}) of $X$ is a maximal anticonnected subset of $X$.

A \emph{path} $P$ of $T$ is a sequence of distinct vertices $p_1,\cdots,p_k$ such that either $k=1$, or for $i,j\in \{1,\cdots,k\}$, $p_i$ is adjacent to $p_j$ if $|i-j|=1$ and $p_i$ is antiadjacent to $p_j$ if $|i-j|>1$. We say that $P$ is a path \emph{from} $p_1$ \emph{to} $p_k$, and the \emph{endpoints} of $P$ are $p_1$ and $p_k$. Under these conditions, let $V(P)=\{p_1,\cdots,p_k \}$, the $\emph{interior}$ of $P$, denoted by $P^*$, is the induced subtrigraph of $P$ with $V(P^*)=V(P) \setminus \{p_1,p_k\}$, and the \emph{length} of $P$ is $k-1$. We say $P$ is \emph{even} (\emph{odd}) if it has even (odd) length. Two paths $P_1$ and $P_2$ are \emph{disjoint} if $V(P_1)\cap V(P_2)=\emptyset$, and they are \emph{internally disjoint} if $V(P_1^*) \cap V(P_2^*)=\emptyset$; $P_1$ is a \emph{subpath} of $P_2$ if $P_1$ is a connected induced subtrigraph of $P_2$. Sometimes we denote $P$ by $p_1 \c p_k$. Notice that, as a graph is also a trigraph, our definition of a path of a graph here is equivalent to a \emph{chordless path} of a graph in some literature. 

A \emph{cycle} in a trigraph $T$ is an induced subtrigraph $H$ of $T$ with vertices $h_1,\cdots, h_k$ such that $k\geq 3$, and for $i,j \in \{1,\cdots,k\}$, $h_i$ is adjacent to $h_j$ if $|i-j|=1$ or $|i-j|=k-1$; a \emph{hole} is a cycle that further satisfies that $h_i$ is antiadjacent to $h_j$ if $1<|i-j|<k-1$. The \emph{length} of a hole (cycle) is the number of vertices in it. Sometimes we denote $H$ by $h_1 \c h_k \d h_1$. An \emph{antipath} (\emph{antihole}) in $T$ is an induced subtrigraph of $T$ whose complement is a path (hole) in $\overline{T}$. 

A \emph{prism} in a trigraph $T$ is an induced subtrigraph $H$ such that the full realization of $H$ is a prism. A trigraph $T$ is \emph{Berge} if it contains no odd hole and no odd antihole. By this definition, $T$ is Berge if and only if $\overline{T}$ is Berge. Also, $T$ is Berge if and only if every realization (semirealization) of $T$ is Berge. An \emph{even pair} in $T$ is a strongly nonadjacent pair $uv \in \binom{V(T)}{2}$ such that every path from $u$ to $v$ in $T$ is even.

\subsection{Basic Trigraphs} \label{sectbasicdef}

Here, we define the classes of basic trigraphs. A trigraph $T$ is \emph{bipartite} if its vertex set can be partitioned into two strongly stable sets, called a {\em bipartition}. A trigraph $T$ is a \emph{line trigraph} if its full realization is the line graph of a bipartite graph and every clique of size at least $3$ in $T$ is a strong clique. A trigraph is a \emph{doubled graph} if it has a \emph{good partition}. A good partition is a partition $(X,Y)$ of $V(T)$  satisfying the following:
\begin{itemize}
    \item Every component of $T\mid X$ has at most two vertices, and every anticomponent of $T \mid Y$ has at most two vertices.
    \item No switchable pair of $T$ is between $X$ and $Y$. 
    \item For every component $C_x$ of $T|X$ and every anticomponent $C_y$ of $T|Y$, every vertex $v$ of $C_x \cup C_y$ is incident with at most one strong edge and at most one strong antiedge between $C_x$ and $C_y$.
\end{itemize}

A trigraph is \emph{basic} if it is either a bipartite trigraph, the complement of a bipartite trigraph, a line trigraph, the complement of a line trigraph, or a doubled trigraph. The following is Theorem 2.3 from \cite{algoTri}:

\begin{theorem}[\cite{algoTri}]\label{lm:close}
Basic trigraphs are Berge, and are closed under taking induced subtrigraphs, semirealizations, realizations, and complementation.
\end{theorem}

\subsection{Decompositions}
We now describe the decompositions for trigraphs. First, a {\em $2$-join} in a trigraph $T$ is a
partition $(X_1, X_2)$ of $V(T)$ such that there exist disjoint sets
$A_1, B_1, C_1, A_2, B_2, C_2 \subseteq V(T)$ satisfying:

\begin{itemize}
\item $X_1=A_1\cup B_1\cup C_1$ and $X_2=A_2\cup B_2\cup C_2$;
\item $A_1, A_2, B_1$ and $B_2$ are non-empty;
\item no switchable pair is between $X_1$ and $X_2$;
\item every vertex of $A_1$ is strongly adjacent to every vertex of
  $A_2$, and every vertex of $B_1$ is strongly adjacent to every
  vertex of $B_2$;
\item there are no other strong edges between $X_1$ and $X_2$; 
\item for $i=1,2$ $|X_i| \geq 3$; and 
\item for $i = 1,2$, if $|A_i| = |B_i| = 1$, then the full realization of
$T|X_i$ is not a path of length two containing the members of $A_i$ and $B_i$.
\end{itemize}

Under these conditions, we say that $(A_1, B_1, C_1, A_2, B_2, C_2)$
is a \emph{split} of $(X_1, X_2)$. A $2$-join is \emph{proper} if
for $i = 1,2$, every component of $T|X_i$ meets both $A_i$ and
$B_i$. A \emph{complement $2$-join} of a
trigraph $T$ is a $2$-join of $\overline{T}$. We need the following fact about $2$-joins (Theorem 2.4 of \cite{algoTri}):

\begin{theorem}[\cite{algoTri}]\label{l:par2Join} 
Let $T$ be a Berge trigraph and $(A_1, B_1, C_1,
A_2, B_2, C_2)$ a split of a proper $2$-join of $T$. Then all paths
with one end in $A_i$, one end in $B_i$ and interior in $C_i$, for
$i=1, 2$, have lengths of the same parity.
\end{theorem}

Next, a partition $(A,B)$ of $V(T)$ is a \emph{skew-partition} if $A$ is not connected and $B$ is not anticonnected. A skew-partition $(A,B)$ is \emph{balanced} if there is no odd path of length greater than one with ends in $B$ and interior in $A$, and there is no odd antipath of length greather than one with ends in $A$ and interior in $B$. Given a balanced skew-partition $(A,B)$, the 4-tuple $(A_1,A_2,B_1,B_2)$ is a
\emph{split of $(A,B)$} if $A_1, A_2, B_1$, and $B_2$ are disjoint
non-empty sets, $A_1 \cup A_2 =A$, $B_1 \cup B_2=B$, $A_1$ is strongly
anticomplete to $A_2$, and $B_1$ is strongly complete to $B_2$. Note
that there exists at least one
split for every balanced skew-partition.

When $(A, B)$ is a skew-partition of a trigraph $T$, we say that $B$ is
a \emph{star cutset} of $T$ if at least one anticomponent of $B$ has size one.
The following is Theorem~5.9 from~\cite{thesis}.

\begin{theorem}[\cite{thesis}]
 \label{starcutset}
 If a Berge trigraph admits a star cutset, then it admits a balanced skew-partition.
\end{theorem}

We will often use the following corollary: 
\begin{theorem}[\cite{thesis}]
\label{no-sc}
If $T$ is a Berge trigraph with no balanced skew-partition, then $T$ does not admit a star cutset.
\end{theorem}

\section{Decomposing Trigraphs}\label{decomp}
\subsection{Decomposing Trigraphs from \texorpdfstring{$\mathcal{F}$}{F}.} In order to handle $2$-join partitions and their complements in Section \ref{BofD}, we define a class of trigraphs that will be useful. 

Let $T$ be a trigraph. Denote by $\Sigma(T)$ the graph with vertex set $V(T)$ and edge set $\sigma(T)$ (the switchable pairs of $T$). 
A \emph{switchable component} of $T$ is a connected component of $\Sigma(T)$
with at least two vertices. Let $\mathcal{F}$ be the class of Berge trigraphs $T$ such that the following hold:
\begin{enumerate}
\item $T$ has at most one switchable component, and the switchable component $D$ of $T$ has at most two edges.
\item If $D$ contains exactly one edge $xy$, then $N(x)\cap N(y)=\emptyset$ in the trigraph $T$. In this case, we say it is a $\emph{small}$ switchable component.
\item Next, assume that $D$ has two edges. Let $v \in V(T)$ be the vertex of degree two in $\Sigma(T)$, denote its neighbors by $x$ and~$y$. Then $v$ is strongly anticomplete to $V(T) \setminus \{v, x, y\}$ in $T$, $x$ is strongly antiadjacent to $y$ in $T$, and $N(x)\cap N(y)=\{v\}$ in $T$. In this case, we say that the switchable component is \emph{light}.
\end{enumerate}
Let $T \in \mathcal{F}$. Then $T$ has at most one switchable component.
If $T$ has exactly one switchable component,  we denote the  vertex set of such a component by $D$. If $T$ has no switchable components
(that is, $\sigma(T)=\emptyset$), we set $D=\emptyset$. In
both cases we refer to $D$ as ``the switchable component of $T$''.
We also use  $D$ to refer to  $\Sigma(T)|D$.
Our class $\mathcal{F}$ of trigraphs is a subclass of the class of the same name studied in \cite{algoTri}, and we make use of several of their results.

\begin{theorem}[\cite{algoTri}]\label{thm:2jnobsp}
Let $T$ be a trigraph from $\mathcal F$ with no balanced
  skew-partition, and let $(A_1,B_1,C_1,A_2,B_2,C_2)$ be a split of a
  $2$-join $(X_1,X_2)$ in $T$. Then the following hold:
 \begin{enumerate}[(i)]
 \item $(X_1,X_2)$ is a proper $2$-join;
 \item if $C_i=\emptyset$, then $|A_i|\geq 2$ and $|B_i|\geq 2$, $i=1,2$;
 \item  $|X_i| \geq 4$, $i=1,2$.
 \end{enumerate}
\end{theorem}

\begin{theorem}[\cite{algoTri}]\label{lemma:t-2join}
Every trigraph in $\mathcal{F}$ is either basic, or admits a proper $2$-join, or admits a proper $2$-join in the complement.
\end{theorem}

\subsection{Favorable Trigraphs}
Let $T$ be a trigraph in $\mathcal{F}$. We say a pair $uv$ of vertices of $T$ is \emph{disjoint from its switchable component} if $D$ is the switchable component of $T$ and $V(D)\cap \{u,v\}$ is empty. In particular, if the switchable component $D$ of $T$ is empty, every pair of vertices is disjoint from its switchable component. A trigraph $T\in \mathcal{F}$ is \emph{favorable} if it satisfies the following conditions:
\begin{enumerate}
    \item $|V(T)|\geq 5$;
    \item $T$ has at least one pair of strongly nonadjacent vertices $uv$ disjoint from $D$; and
    \item if $D$ is small and $V(D)=\{x,y\}$, then at least one of $T\setminus (D\cup N(x))$ or $T\setminus (D\cup N(y))$ is not a clique.
\end{enumerate}
A trigraph is $\emph{unfavorable}$ if it is not favorable. By this definition, if $T$ is complete, then $T$ is unfavorable; if $T$ is a graph with at least five vertices and is not complete, then $T$ if favorable as it has empty switchable component. Notice that condition $(2)$ and $(3)$ of being a favorable trigraph are also necessary conditions for trigraphs to have even pairs disjoint from the switchable component. 

Next, we will show that, with a few exceptions, a trigraph $T$ in $\mathcal{F}$ with no balanced skew-partition and no antihole is favorable. Further, we prove in section \ref{Basic} that a basic favorable trigraph has an even pair disjoint from its switchable component. Both results are essential for handling $2$-joins in section \ref{BofD}.

\begin{theorem}\label{switcompdom}
Let $T$ be a trigraph in $\mathcal{F}$ with no balanced skew-partition and no antihole of length six. If $T$ is unfavorable, then either $T$ is complete or $|V(T)|\leq 5$.
\end{theorem}

\Proof 
We may assume that $|V(T)|>5$. Let $D$ be the switchable component of $T$, and let $T'=T\setminus V(D)$ be the induced subtrigraph of $T$. If $D$ is small, we denote the pair by $x$ and $y$; if $D$ is light, we denote the vertex of degree two in $\Sigma(T)$ by $v$ and its neighbors by $x$ and $y$. Therefore, we can partition $V(T')$ into three sets: $T_1=T' 
\setminus (N(x) \cup N(y))$, $T_2=T'|N(x)$, $T_3=T'|N(y)$.

First, suppose that $D$ is a light switchable component. Since $D$ is unfavorable, it follows that $V(T) \setminus V(D)$ is a clique. If both $T_2$ and $T_3$ are nonempty, then $x\d a \d b \d y \d v \d x$ with $a\in T_2$ and $b\in T_3$ is a hole of length five, contradicting that $T$ is Berge, so we may assume up to symmetry that $T_3 = \emptyset$. Since $|V(T)| \geq 6$, it follows that $T_1 \cup T_2$ contains three distinct vertices including $s$ and $t$. If $T_2=\emptyset$, then $\{s,t\}$ is a star cutset, contradicting $\ref{no-sc}$. So we may assume $t\in T_2$. Now, $V(T) \setminus \{v,y,s\}$ is a star cutset, again contradicting Theorem~\ref{no-sc}.

Therefore, $D$ is a small switchable component. Now, if $T_1\neq \emptyset$, $T_2$ is strongly complete to $T_3$ since $T$ contains no hole of length five. In this case, $T_2\cup T_3$ is a star cutset, contradicting Theorem~\ref{no-sc}. Thus, $T_1 =\emptyset$. By the definition of unfavorable, since $T$ is not complete, it follows that both $T_2$ and $T_3$ are cliques. As $T$ has at least five vertices, at least one of $T_2$ or $T_3$ has more than two vertices. Without loss of generality, suppose $T_2$ contains two distinct vertices. Let $s$ be the vertex in  $T_2$ such that $|N(s) \cap T_3|$ is the maximum, and let $t$ be a vertex in $T_2$ distinct from $s$. By Theorem~\ref{no-sc}, we may assume $N(s)\cup \{s\} \setminus \{t\}$ is not a star cutset. It follows that there exists a vertex $p\in T_3 \setminus N(s)$ adjacent to $t$. By maximality of $|N(s) \cap T_3|$, there exists $q\in N(s)\cap T_3$ such that $q$ is not connected to $t$. Now, $T| \{x,y,s,t,p,q\}$ is an antihole of length of six, a contradiction. This completes the proof.
\bbox

\section{Even Pairs in Basic Trigraphs} \label{Basic}
The goal of this section is to prove the following theorem by analyzing each class of basic trigraph:

\begin{theorem}\label{thm:basicdisj}
Let $T$ be a basic trigraph in $\mathcal{F}$ with no odd prism and no antihole. Then the following statements hold:
\begin{enumerate}
    \item $T$ is either complete or has an even pair.
    \item If $T$ is favorable, then $T$ has an even pair disjoint from its switchable component. 
\end{enumerate}
\end{theorem}

\subsection{Bipartite Trigraph}
Let $T$ be a bipartite trigraph with bipartition $(A, B)$, where $A$ and $B$ are strongly stable sets. We have the following observation.

\begin{theorem}\label{lm:bip}
Let $T$ be a bipartite trigraph in $\mathcal{F}$. Then the following statements hold:
\begin{enumerate}
    \item $T$ is either complete or has an even pair.
    \item If $T$ is favorable, then $T$ has an even pair disjoint from its switchable component. 
\end{enumerate}
\end{theorem}

\begin{Proof} By the definition of bipartite trigraph, it holds that $T$ is complete or has an even pair, so the first statement holds. For the second statement, suppose that $T$ is favorable and has a nonempty switchable component $D$. If either $A'=A \setminus V(D)$ or $B'=B\setminus V(D)$ contains at least two vertices, then any two vertices $a_1,a_2\in A'$ (or $b_1,b_2\in B'$) form an even pair disjoint from the switchable component, so we may assume that $|A'| = |B'| = 1$. Since $T$ is favorable, it follows that $|V(T)| \geq 5$. Thus $|V(D)| \geq 3$, and so $T$ has a light switchable component, and $|V(T)|=5$. Assume up to symmetry that $A=\{v,a\}$ and $B=\{x,y,b \}$, where $v$ is the vertex of degree two in $D$, and $x$ and $y$ are neighbors of $v$ in $D$. Since $T$ is favorable, it follows that $ab$ is a strong antiedge and $b$ is strongly anticomplete to $V(T) \setminus \{b\}$. Now, $a$ and $b$ are in disjoint connected components of $T$, so $ab$ is an even pair of $T$ disjoint from the switchable component. \bbox
\end{Proof}

\subsection{Line Trigraph}
Let $T$ be a line trigraph, and let $H$ be the bipartite graph such that its line graph, denoted by $L(H)$, is the full realization of $T$. Let $(A,B)$ be a bipartition of $H$. A pair $(a_1b_1,a_2b_2)$ of disjoint edges in $H$ with $a_1,a_2\in A $ and $b_1,b_2\in B$ is a \emph{good pair} of $H$ if both of the followings are satisfied:
\begin{itemize}
    \item Every path $P_1$ with endpoints $a_1$ and $a_2$ satisfies $V(P_1)\cap \{b_1,b_2\}\neq \emptyset$; and
    \item every path $P_2$ with endpoints $b_1$ and $b_2$ satisfies $V(P_2)\cap \{a_1,a_2\}\neq \emptyset$.
\end{itemize}
We prove that a good pair in $H$ corresponds to an even pair in $T$. This is analogous to a result by Hougardy in \cite{HOUGARDY199517}.

\begin{prop}\label{lm:goodpair}
Let $H$ be a bipartite graph, let $(a_1b_1,a_2b_2)$ be a good pair of $H$, and let $u$ and $v$ be the vertices in $L(H)$ that represent $a_1b_1$ and $a_2b_2$, respectively. Let $T$ be a trigraph such that $L(H)$ is the full realization of $T$. Then, $uv$ is an even pair in $T$.
\end{prop}
\Proof First, note that $uv$ is a strong antiedge in $T$, as $a_1b_1$ and $a_2b_2$ are disjoint in $H$. Suppose that there is an odd path $P$ from $u$ to $v$ in $T$. Then, $P$ corresponds to an even path $Q$ in $H$ such that $L(Q)=P$. Therefore, up to symmetry, we may assume that $Q=a_1 \d b_1 \d \cdots \d b_2 \d a_2$, and let $Q'$ be the subpath of $Q$ with endpoints $b_1$ and $b_2$. Thus, $V(Q')\cap \{a_1,a_2\}=\emptyset$. However, this contradicts that $(a_1b_1,a_2b_2)$ is a good pair. \bbox

Next, we show that forbidding odd prisms guarantees even pairs in line trigraphs. Let $H$ be a bipartite graph. (Note that the following theorems consider all subgraphs of $H$, which are not necessarily induced subgraphs.) A (not necessarily induced) path $Q$ of $H$ is a \emph{chord path} if its endpoints are contained in the vertex set of a cycle $C$ in $H$ and $V(Q^*)\cap V(C)=\emptyset$. A \emph{path along the cycle} $C$ is an induced subgraph of $C$ that is a path in $H|V(C)$. A graph is \emph{series-parallel} if and only if it has no subgraph isomorphic to a $K_4$-minor. An \emph{even theta} is a graph composed of three internally disjoint even paths with the same endpoints. 

\begin{prop}\label{lm:k4}
Let $T$ be a line trigraph with no odd prism, and let $H$ be a bipartite graph such that $L(H)$ is the full realization of $T$. Then, $H$ has no subgraph isomorphic to an even theta, and $H$ is series-parallel. 
\end{prop}

\Proof Since the line graph of an even theta is an odd prism, it follows that $H$ contains no even theta as a subgraph. As $H$ is bipartite, all cycles in $H$ have even length. It follows that a chord path $P$ of a cycle $C$ in $H$ must has odd length, and the endpoints of $P$ in $C$ divide $C$ into two odd paths along the cycle. Suppose the contrary that $H$ is not series-parallel and thus has a subgraph isomorphic to a $K_4$-minor. Since $K_4$ has maximum degree three, it follows that $H$ has a subgraph $J$ isomorphic to a $K_4$-subdivision. Let $a,b,c,d$ be the vertices of degree three of $J$, and let $P_1$, $P_2$, $P_3$, $P_4$, $P_5$, $P_6$ denote the paths with endpoints $(a,b)$, $(b,c)$, $(c,d)$, $(d,a)$, $(b,d)$, and $(a,c)$, respectively, in $J$. Notice that each $P_i$ is a chord path, so they are all odd. Now, $P_1\cup P_4 \cup P_5$ is an odd cycle, contradicting that $H$ is bipartite. \bbox

Finally, we prove the main result of this subsection. 

\begin{theorem}\label{lm:linetri}
Let $T$ be a line trigraph in $\mathcal{F}$ with no odd prism. The following statements hold:
\begin{enumerate}
    \item $T$ is either complete or has an even pair.
    \item If $T$ is favorable, then $T$ has an even pair disjoint from its switchable component.
\end{enumerate}
\end{theorem}

\Proof Let $H$ be a bipartite graph such that $L(H)$ is the full realization of $T$. We may assume that $H$ is connected and that $T$ is not complete. If $T$ has a nonempty switchable component $D$, let $J$ be the subgraph of $H$ such that $T|V(L(J))=D$. In particular, $J$ is a path $p_1 \d \cdots \d p_k$ of length either two or three. Thus, following the notation for a path, we call $p_1$ and $p_k$ the endpoints of $J$ and denote $V(J)\setminus \{p_1,p_k\}$ by $V(J^*)$. Also, note that any vertex $v\in V(J^*) $ has degree at most two in $H$: Otherwise, the line graph induced by the edges adjacent to $v$ is a clique $K$ of size at least three, and $T|V(K)$ contains a switchable pair, which contradicts the definition of a line trigraph. 

By Theorem~\ref{lm:goodpair}, to prove the first statement, it suffices to find a good pair in $H$. Also, to prove the second statement, it suffices to find a good pair in $H\setminus V(J^*)$. Thus, in the following discussion, the proof is completed when the corresponding good pair is found.

\textbf{Case 1: $H$ is a tree.} Since $T$ is not complete, it follows that $H$ is not a star, so $H$ has a path $a_1\d b_1 \d a_2 \d b_2$ of length three. Now, $(a_1b_1,a_2b_2)$ is a good pair. This proves the first statement for this case. 

Next, suppose that $T$ is favorable and has a nonempty switchable component $D$. Let $x$ and $y$ be the endpoints of $J$, and let $H_x$ and $H_y$ be the components of $H\setminus V(J^*)$ containing $x$ and $y$ correspondingly. It suffices to show that $H_x\cup H_y$ contains a good pair. If either $H_x$ or $H_y$ contains a path $a_i \d b_i \d a_j \d b_j$ of length three, then $(a_ib_i,a_jb_j)$ is a good pair. Thus, we may assume both $H_x$ and $H_y$ are either empty or isomorphic to a star. If $D$ is small, then $T$ contradicts the third condition of being favorable. So we may assume $D$ is light. By the second condition of being favorable, both $H_x$ and $H_y$ are nonempty. In particular, $H_x$ contains an edge $xx'$, and $H_y$ contains an edge $yy'$. Now, $(xx',yy')$ is a good pair. This completes the proof of the second statement for this case.

\textbf{Case 2: $H$ has a cycle of length at least six}. Let $C=a_1\d b_1 \d \cdots \d a_k \d b_k \d a_1$ where $k\geq 6$ be a cycle (not necessarily induced) of maximum length in $H$. If $C$ has no chord path, then every pair of disjoint edges $(a_ib_i,a_jb_j)$ is a good pair. In particular, as $|E(J)|\leq 3$, there is a good pair in $C\setminus V(J^*)$. Thus, we may assume that $C$ has a chord path $P$. By Theorem~\ref{lm:k4}, $P$ is odd and has ends $a_i$ and $b_j$ for $1 \leq i \leq j \leq k$. Let $Q_1$ and $Q_2$ be the two disjoint paths along the cycle $C$ with endpoints $a_i$ and $b_j$. We may assume by symmetry that $E(J)\cap E(Q_1)=\emptyset$ as any $v\in V(J^*)$ has degree two. Now, to prove both statements for this case, it suffices to show that $Q_1$ contains a good pair.

Let $S_1$ be a minimal subpath of $Q_1$ such that the endpoints of $S_1$ are joined by a chord path of $C$, and let this chord path be $P'$. Thus, $S_1$ has odd length. If $S_1$ has length one, then $P'\cup (C \setminus S_1)$ is a longer cycle, a contradiction. So $S_1=a_t \d b_t \d \cdots \d a_s\d b_s$ has length at least three. Further, if there is a chord path $P''$ of $C$ with exactly one endpoint in $V(S_1^*)$, then $C\cup P' \cup P''$ forms a $K_4$ minor, contradicting Theorem~\ref{lm:k4}. Therefore, there is no path in $H\setminus \{a_s,a_t\}$ with endpoints $b_s$ and $b_t$, and there is no path in $H\setminus \{b_s,b_t\}$ with endpoints $a_s$ and $a_t$. So $(a_tb_t,a_sb_s)$ is a good pair of $H$ contained in $Q_1$. This completes the proof.

\textbf{Case 3: All the cycles in $H$ have length four.} Let $C=a_1 \d b_1 \d a_2 \d b_2 \d a_1$ be a cycle of length four in $H$. By Theorem~\ref{lm:k4}, there is no chord path of $C$ with endpoints $a_1$ and $a_2$ (or $b_1$ and $b_2$). Also, if there is a chord path with endpoints $a_i$ and $b_j$ with $i,j\in \{1,2\}$, then $G$ contains a cycle of length greater than four, a contradiction. Thus, there is no path in $H\setminus \{a_1,a_2\}$ with endpoints $b_1$ and $b_2$, and there is no path in $H\setminus \{b_1,b_2\}$ with endpoints $a_1$ and $a_2$. So $(a_1b_1, a_2b_2)$ is a good pair of $H$. In particular, this proves that in this case every cycle $C$ in $H$ contains a good pair, and thus the first statement follows. 

Now, suppose $T$ is favorable with nonempty switchable component $D$. We may assume that $E(J)\cap E(C)\neq \emptyset$, and $H\setminus V(J^*)$ contains no cycle. As any vertex $v\in V(J^*)$ has degree two, we have $E(J)\subseteq E(C)$ and $V(J)\subseteq V(C)$. Thus, we may assume that $H\setminus V(C)$ is a tree. If $D=\{x,y,v\}$ (where $v$ has degree two in the switchable component) is a light switchable component, then $N(x) \cap N(y)\neq \emptyset$ as the endpoints of $J$ are adjacent, contrary to the fact that $T\in \mathcal{F}$. Therefore, we may assume that $D$ is small and $J= a_1\d b_1 \d a_2$. As $T$ is favorable, $T\setminus V(D)$ is not a clique, which means that there is an edge $a_tb_t$ in $H\setminus V(J^*)$ such that $b_t\neq b_2$. Also, since $b_1\in V(J^*)$ has degree two in $H$, we have $b_t\neq b_1$. If $a_t=a_1$, then $\{a_tb_t,a_2b_2\}$ is a good pair in $H\setminus V(J^*)$. Thus, by symmetry, we may assume that $\{a_t,b_t\}\cap V(C)=\emptyset$, and every edge between $C$ and $H\setminus V(C)$ has $b_2$ as a vertex. In this case, as $C$ has no chord path and $H\setminus V(C)$ is a tree, $\{a_tb_t, a_2b_2\}$ is a good pair in $H\setminus V(J^*)$. This completes the proof. \bbox

\subsection{Complement of a Bipartite Trigraph and Complement of a Line Trigraph} A \emph{diamond} in a trigraph $T$ is an induced subtrigraph $H$ such that the full realization of $H$ is $K_4$ minus an edge. A \emph{claw} in a trigraph $T$ is an induced subtrigraph $H$ such that the full realization of $H$ is the complete biparite graph $K_{1,3}$. We will need the following characterization of line trigraph, which is a generalization of the main theorem of \cite{harary1974line}.

\begin{prop}\label{lm:nocd}
Let $T$ be a line trigraph. Then, $T$ has no induced subtrigraph isomorphic to a diamond or a claw.
\end{prop}
\Proof
Suppose that $T$ contains a diamond or a claw, then the full realization of $T$ contains a diamond or a claw as an induced subgraph. By definition, the full realization of $T$ is a line graph. This contradicts to the main theorem of \cite{harary1974line}, which states that a line graph of a bipartite graph is (claw,diamond)-free.
\bbox

Basic trigraphs which are the complement of a bipartite trigraph and the complement of a line trigraph share the following key property. 

\begin{prop}\label{lm:3path}
Let $T$ be the complement of a bipartite trigraph or the complement of a line trigraph. Then, a path $P$ of odd length in $T$ has length at most three.
\end{prop}

\Proof First, if $T$ is the complement of a bipartite trigraph, then for all $X\subseteq V(T)$ with $|X|\geq 3$, there exists an edge with both ends in $X$. Therefore, the path of maximal length in $T$ has length three, so the result follows. Now, we may suppose that $T$ is the complement of a line trigraph, and $P$ is a path of $T$ of length at least five. In this case, $\overline{T|V(P)}$ contains a diamond, which contradicts Theorem~\ref{lm:nocd}. This completes the proof. \bbox

\begin{prop}\label{lm:swcomp}
Let $T$ be a trigraph in $\mathcal{F}$ such that $T$ is either the complement of a bipartite trigraph or the complement of a line trigraph. If $T$ is favorable, then either $T$ is a graph, or $T$ has an even pair disjoint from the switchable component.
\end{prop}

\Proof Recall that by the definitions, every clique of size at least three is a strong clique in line trigraphs and bipartite trigraphs. If $T$ has a light switchable component $D=\{x,y,v\}$, then, $\overline{T|D}$ is a clique of size three with two switchable pairs, contradicting that $T$ is the complement of a bipartite trigraph or the complement of a line trigraph. Thus, we may assume that $T$ has a small switchable component $D=\{x,y\}$. Suppose that $T$ is the complement of a bipartite trigraph with bipartition $(A,B)$. By definition, $T|A$ and $T|B$ are strong cliques, so we may assume $x\in A$ and $y\in B$ up to symmetry. In this case, $T\setminus (D \cup N(x))\subseteq B$ and $T\setminus (D \cup N(y))\subseteq A$ are both cliques, contradicting that $T$ is favorable. 

Now, we may assume that $T$ is the complement of a line trigraph with a small switchable component. If there is a vertex $v$ contained in $T\setminus (N(x)\cup N(y))$, then $\overline{T|\{x,y,v\}}$ is a clique of size three but not a strong clique, contradicting the definition of line trigraph. So $T\setminus (N(x)\cup N(y))=\emptyset$. If there are two vertices $s,t\in N(x)$ (or $s,t\in N(y)$) such that $st$ is an edge in $T$, then $\overline{T|\{x,s,t,y\}}$ is a claw, contradicting Theorem~\ref{lm:nocd}. Therefore, $T|N(x)$ and $T|N(y)$ are stable sets. Since $T$ is favorable, we may assume up to symmetry that $|N(x)|\geq 2$. Let $s$ and $t$ be two vertices in $N(x)$, and we claim that $\{s,t\}$ is an even pair: Suppose not, then there is an odd path $s\d v_1 \d v_2 \d t$ of length three by Theorem~\ref{lm:3path}. As $T|N(x)$ is a stable set, $\{v_1,v_2\}\subseteq T\setminus (N(x)\cup \{x,y\})=N(y)$. So $v_1v_2$ is an edge in $T|N(y)$, contradicting that $T|N(y)$ is a stable set. This completes the proof. \bbox

The proof of following proposition is inspired by the main idea of \cite{MAFFRAY}.
\begin{theorem} \label{lm:unifed}
Let $T$ be a trigraph in $\mathcal{F}$ with no antihole such that $T$ is either the complement of a bipartite trigraph or the complement of a line trigraph. The following statements hold:
\begin{enumerate}
    \item $T$ is either complete or has an even pair.
    \item If $T$ is favorable, then $T$ has an even pair disjoint from its switchable component.
\end{enumerate}
\end{theorem}

\Proof By Theorem~\ref{lm:swcomp}, it suffices to prove the first statement. Let $T$ be the vertex-minimal counterexample. We may assume that $T$ is not complete. Let $M$ be a maximal anticonnected set in $T$ such that there are at least two nonadjacent vertices in $V(T) \setminus M$ that are complete to $M$. Notice that $M$ is nonempty: since $T$ is not complete, it holds that $T$ contains at least one path of length at least two. Let $C(M)$ be the set of all vertices that are complete to $M$. By Theorem~\ref{lm:close}, each class of basic trigraphs is closed under taking induced subtrigraphs. Thus, since $T$ is minimal, it follows that $C(M)$ has an even pair $\{a,b\}$ as $C(M)$ is not complete by our construction. 

Suppose the contrary that $\{a,b\}$ is not an even pair in $T$. Thus, by Theorem~\ref{lm:3path}, there is a path $P=a\d c \d d \d b$ of length three in $T$. Since $\{a, b\}$ is complete to $M$, it follows that $V(P)\cap M=\emptyset$. First, suppose $\{c, d\} \subseteq V(T)\setminus (M \cup C(M))$. Since both $c$ and $d$ are not in $C(M)$, it follows that $c$ and $d$ each has at least one strong antineighbor in $M$. So there exists an antipath $Q$ with ends $c$ and $d$ and $Q^*\in M$. Then, $c \d Q \d d \d a \d b \d c$ is an antihole of length at least five in $T$, a contradiction. Thus, we may assume up to symmetry that $c \in C(M)$. Since $\{a, b\}$ is an even pair in $C(M)$, it follows that $P \not \subseteq C(M)$, and so $d \in V(T)\setminus (M \cup C(M))$. Now, $M \cup \{d\}$ is also an anticonnected set in $T$, and $\{c,b\}$ is a pair of nonadjacent vertices complete to $M \cup \{d\}$. This contradicts that $M$ is maximal. Therefore, $\{a,b\}$ is an even pair in $T$. This completes the proof. \bbox

\subsection{Doubled Graph} We first state a proposition regarding even pairs in doubled graphs.

\begin{prop}\label{dgprop}
Let $T$ be a doubled graph with good partition $(X,Y)$. Then the following two statements hold:
\begin{enumerate}
    \item Let $C_y=\{y\}$ be an anticomponent of size one in $Y$. If $y$ has an antineighbor $x\in X$, then $xy$ is an even pair. In particular, if $X$ has an edge $x_1x_2$, then $T$ has an even pair.
    \item If $T|Y$ has a strong antiedge $y_1y_2$ and $T|X$ has no edge, then $y_1y_2$ is an even pair.
\end{enumerate}
\end{prop}

\Proof
\begin{enumerate}
    \item In this case, $y$ is complete to $N(x)$, so all paths between $x$ and $y$ have length $2$. If $X$ contains an edge $x_1x_2$, then either $x_1$ or $x_2$ is an antineighbor of $y$, so one of $x_1y$ or $x_2y$ is an even pair.
    \item In this case, $N(y_1)\cap X$ and $N(y_2)\cap X$ are two disjoint stable sets that partition $X$. So there is no path from $y_1$ to $y_2$ whose interior is contained in $X$. Since $\{y_1,y_2\}$ is complete to other vertices in $Y$, all paths from $y_1$ to $y_2$ have length $2$. \bbox
   
\end{enumerate} 

\begin{theorem}\label{thm:dbg}
Let $T$ be a doubled graph in $\mathcal{F}$ with no antihole of length six. The folowing statements hold:
\begin{enumerate}
    \item $T$ is either complete, or has an even pair. 
    \item If $T$ is favorable, then $T$ has an even pair disjoint from its switchable component.
\end{enumerate}
\end{theorem}
\Proof We may assume that $T$ is not complete and is connected, as any two vertices taken from distinct connected components form an even pair. Let $(X, Y)$ be a good partition of $T$. Note that by the definition of doubled graph, every switchable component of $T$ is either an edge of $T|X$ or an edge of $T|Y$. If $T$ has a switchable component $D$, it must be small: otherwise, $T|D$ is both a component and an anticomponent of size $3$, contradicting that $T$ is a doubled graph.

\textbf{Case 1: $T|X$ has a component $C_1 = \{x_1, x_2\}$ of size two.} If $T|Y$ is empty, then $T = C_1$ and $T$ is a clique, so we may assume that $T|Y$ is nonempty. If $T|Y$ has two distinct anticomponents $C_2,C_3$ of size two, $C_1\cup C_2 \cup C_3$ is an antihole of length six, a contradiction. Therefore, $T|Y$ has at most one anticomponent of size two. 

First, suppose $T|Y$ has an anticomponent $C_4=\{v\}$ of size one. By symmetry, we may assume that $v$ is strongly adjacent to $x_1$ and strongly antiadjacent to $x_2$. By $(1)$ of Theorem~\ref{dgprop}, it follows that $\{v, x_2\}$ is an even pair of $T$, and this proves the first statement for this subcase. Next, assume that $T$ has a small switchable component $D$. If $\{x_1, x_2\}$ is not the switchable component of $T$, then $\{v, x_2\}$ is an even pair disjoint from its switchable component, so assume that $\{x_1, x_2\}$ is the switchable component. By the structure of a doubled graph, $Y$ is partitioned by $Y_1=N(x_1) \setminus \{x_2\}$ and $Y_2=N(x_2)\setminus \{x_1\}$, and both $Y_1$ and $Y_2$ are cliques. Since $T$ is favorable, there must exists a vertex $x_3$ in  $T|X \setminus \{x_1, x_2\}$ such that $x_3$ has an antineighbor $y_1 \in T|Y$. If $\{y_1\}$ is an anticomponent of size one in $T|Y$, $\{x_3,y_1\}$ is an even pair disjoint from the switchable component by $(1)$ of Theorem~\ref{dgprop}. Therefore, we may assume that $y_1$ is in the anticomponent $C_5=\{y_1,y_2\}$ of size two in $T|Y$. As $C_5$ is the only anticomponent of size two in $Y$, we may also assume that $X\setminus \{x_1,x_2\}$ is complete to $Y\setminus \{y_1,y_2\}$. By the structure of a doubled graph, $X\setminus \{x_1,x_2\}$ is a nonempty stable set. Now, $\{x_3,y_1\}$ is an even pair disjoint from the switchable component $\{x_1,x_2\}$: $y_1$ is complete to $N(x)\setminus \{y_2\}$, so a path from $x_3$ to $y_1$ either has length two, or is exactly $x_3\d y_2\d x_2\d x_1\d y_1$, which has length four. This proves the second statement for this subcase. 

Therefore, we may assume that $T|Y$ contains an antiedge $y_1y_2$ and $Y=\{y_1,y_2\}$. By symmetry, we may assume that $x_i$ is strongly adjacent to $y_i$ for $i=1,2$. Notice that all paths from $y_1$ to $y_2$ have length three. Every path $P$ from $y_1$ to $x_2$ goes through either $x_1$ or $y_2$. If $x_1\in P$, then $P=y_1\d x_1\d x_2$ has length two; if $y_2\in P$, then $P=y_1\d \cdots \d y_2 \d x_2$ has length four. Therefore, $y_1x_2$ is an even pair, and $x_1y_2$ is also an even pair by symmetry. Thus, this proves the first statement for this subcase. Next, assume that $T$ is favorable. We may also suppose that either $\{x_1,x_2\}$ or $\{y_1,y_2\}$ is the switchable component. As $|V(T)|\geq 5$, there exists a vertex $x_3\in X\setminus \{x_1,x_2\}$. By symmetry, we may assume $x_3$ is strongly adjacent to $y_1$ and strongly antiadjacent to $y_2$. Suppose $x_3$ has a neighbor $x_4\in X\setminus \{x_1,x_2\}$. Then, $x_3y_2$ and $x_4y_1$ are even pairs by the same argument above. Also, $\{x_1,x_3\}$ and $\{x_2,x_4\}$ are even pairs. So at least one of them is disjoint from the switchable component, which is either $\{x_1, x_2\}$ or $\{y_1,y_2\}$. Thus, we may suppose $x_3$ has no neighbor in $X$. Then, both $\{x_3,y_2\}$ and $\{x_1,x_3\}$ are even pairs: $x_1$ is complete to $N(x_3)=\{y_1\}$, and every path from $x_3$ to $y_2$ that contains $y_1$ has length four. So at least one even pair is disjoint from the switchable component. This completes the proof of both statements for Case 1.

\textbf{Case 2: $T|X$ has no edge.} It follows that the switchable component of $T$ is contained in $T|Y$. By $(2)$ of Theorem~\ref{dgprop}, we may assume that $T|Y$ is a clique. Suppose that there is no switchable pair in $T$. Then $Y$ is a strong clique, and thus $X$ is complete to $Y$. Since $T$ is not complete, there exist nonadjacent vertices $x_1, x_2 \in X$. Now, $\{x_1,x_2\}$ is an even pair of $T$. Next, assume that there exists a switchable pair $y_1y_2$ in $T$. As $T\in \mathcal{F}$, $N(y_1)\cap N(y_2)=\emptyset$, and thus $Y=\{y_1,y_2\}$. Since $T$ is not complete, there exists a vertex $x_1\in X \cap N(y_1)$. Now, $\{x_1, y_2\}$ is an even pair. This proves the first statement for Case $2$. 

Suppose $T$ is favorable with switchable pair $y_1y_2$. Then, there exists $x_3,x_4\in X$ such that either $\{x_3,x_4\}\subseteq N(y_1)$ or $\{x_3,x_4\}\subseteq N(y_2)$. In either cases, $\{x_3, x_4\}$ is an even pair because $N(x_3)=N(x_4)$. This completes the proof. \bbox

\subsection{Proof of Theorem \ref{thm:basicdisj}} By Theorems~\ref{thm:basicdisj}, \ref{lm:linetri}, \ref{lm:unifed}, and \ref{thm:dbg}, we have checked that the statements hold for each class of basic trigraphs. So Theorem~\ref{thm:basicdisj} follows.

\section{Even Pairs in Non-Basic Trigraphs}\label{BofD}
\subsection{Block of Decomposition} To handle $2$-join partitions and their complements, we need the following definitions and a theorem regarding trigraphs with no balanced skew-partition from \cite{algoTri}. 
A set $X\subseteq V(T)$ is a \emph{fragment} of a trigraph $T$ if $(X,V(T)\setminus X)$ is a proper $2$-join of $T$. A proper $2$-join is \emph{even} or \emph{odd} according to the parity of the paths described in Theorem~\ref{l:par2Join}. The \emph{block of decomposition $T_X$} with respect to a fragment $X$ is defined as follows. Let $X_1 = X$, $X_2 = V(T) \setminus X$, and $(A_1, B_1, C_1, A_2, B_2, C_2)$ be a split of the proper 2-join $(X_1, X_2)$. 
\begin{itemize} 
\item {\bf Case 1: $(X_1,X_2)$ is odd}. We build the block of decomposition $T_X=T_{X_1}$ as follows: starting with $T|X_1$, we add two vertices $a$ and $b$, called {\em marker vertices}, such that $a$ is strongly complete to $A_1$, $b$ is strongly complete to $B_1$, $ab$ is a switchable pair, and there are no other edges between $\{a,b\}$ and $X_1$. Note that $\{a,b\}$ is a small switchable component of $T_X$, and we call it the \emph{marker component} of $T_X$.
    
\item {\bf Case 2: $(X_1, X_2)$ is even}. We build the block of decomposition $T_X=T_{X_1}$ as follows: starting with $T|X_1$, we add three vertices $a, b, c$, called \emph{marker vertices}, such that $a$ is strongly complete to $A_1$, $b$ is strongly complete to $B_1$, $ac$ and $bc$ are switchable pairs, and there are no other edges between $\{a,b,c\}$ and $X_1$. Note that $\{a,b,c\}$ is a light switchable component of $T_X$, and we call it the \emph{marker component} of $T_X$.
\end{itemize}

In both cases, we say that $a$ and $b$ are the {\em ends} of the marker component. Again, as our class $\mathcal{F}$ is a subclass of the class of the same name studied in \cite{algoTri}, we make use of the following result.

\begin{theorem}[\cite{algoTri}]\label{thm:frag}
If $X$ is a fragment of a trigraph $T$ in $\mathcal{F}$ with no balanced skew-partition, then the block of decomposition $T_X$ is Berge and has no balanced skew-partition.
\end{theorem}

\begin{theorem}\label{lm:opahinfrag}
Let $X$ be a fragment of a trigraph $T$ in $\mathcal{F}$ with no balanced skew-partition, no odd prism, and no antihole. Then the block of decomposition $T_X$ is Berge, and $T_X$ has no balanced skew-partition, no odd prism, and no antihole of length at least five.
\end{theorem}

\Proof By Theorem~\ref{thm:frag}, it suffices to show $T_X$ has no odd prism and no antihole of length at least six. Let $M$ denote the marker component of $T_X$. Suppose the contrary that $T_X$ has an odd prism $Q$. If $V(Q) \cap M = \emptyset$, then $Q$ is an odd prism in $T$, a contradiction. Therefore, we may assume up to symmetry between the ends of the marker component $a$ and $b$ that $a \in V(Q) \cap M$. Suppose that $N(a)\cap M \not \subseteq V(Q)$.  Let $y \in A_2$ and let $Q' = (Q \setminus \{a\}) \cup \{y\}$. If $V(Q') \subseteq V(T)$, then $Q'$ is an odd prism of $T$, a contradiction. Therefore, $M = \{a, v, b\}$, $a$ is not adjacent to $b$, and $b \in V(Q')$. Let $z \in B_2$ and let $Q'' = (Q' \setminus \{b\}) \cup \{z\}$. Now, $Q''$ is an odd prism of $T$, a contradiction. Therefore, $N(a)\cap M  \subseteq V(Q)$.   As any vertex in $Q$ has degree at least two, and every internal vertex of $M$ has degree exactly two in
$T$, it follows that $M\subseteq V(Q)$, and $M$ is contained in a rung of the prism $Q$. Let $(A_1, B_1, C_1, A_2, B_2, C_2)$ be a split of $(X, V(T) \setminus X)$. Since $(X, V(T) \setminus X)$ is a proper 2-join, it follows that there is a path $P$  of $T$ with ends in $A_2$ and $B_2$ and interior in $C_2$ such that $P$ has the same parity as $M$. Now, $(Q \setminus M) \cup P$ is an odd prism of $T$, a contradiction. Therefore, $T_X$ does not contain an odd prism. 

Next, suppose that $T_X$ contains an antihole, and let $H=v_1 \d \cdots \d v_k \d v_1$ be the shortest antihole in $T_X$. Since $T_X$ is Berge and an antihole of length six is an odd prism, we may assume that $k\geq 7$. When $T_X$ has the marker component $\{ a,b,c\}$, it follows that $c\notin V(H)$ because $c$ is strongly anticomplete to $T_X|X$. If $|V(H)\cap M|=1$, we may assume by symmetry that $a\in V(H)\cap M$. Now, we may replace $a$ by a vertex $a'\in A_2$, and $T|(V(H)\setminus \{a\}) \cup  \{a'\}$ is an antihole of the same length in $T$, a contradiction. Therefore, we may assume that $V(H)\cap M=\{a,b\}$, and let $a=v_i$ and $b=v_j$ where $i<j$. First, suppose $ab$ is an antiedge of the antihole such that $i-j=1$ or $i-j=k-1$. Because $\{a, b\}$ is strongly anticomplete to $C_1$, it follows that no vertex of the antihole is contained in $C_1$. Also, at most one vertex of $H$ is in $A_1$ and at most one vertex of $H$ is in $B_1$, since $a$ is strongly anticomplete to $B_1$ and $b$ is strongly anticomplete to $A_1$. Therefore, $k\leq 4$, a contradiction. Hence, we may suppose that $1<i-j<k-1$. However, as $k\geq 7$, at least one of $v_i \d v_{i+1} \d \cdots \d v_j \d v_i$ or $v_j \d v_{j+1} \d  \cdots \d v_i \d v_j$ is an antihole and has length less than $k$, which is a contradiction that $H$ is the shortest antihole in $T_X$. Therefore, $T_X$ does not contain an antihole. This completes the proof. \bbox

\begin{theorem}\label{thm:epinfrag}
Let $X$ be a fragment of a Berge trigraph $T$. If the block of decomposition $T_X$ has an even pair $uv$ disjoint from its marker component, then $uv$ is also an even pair in $T$.
\end{theorem}

\Proof
Let $(A_1,B_1,C_1,A_2,B_2,C_2)$ be a split of a proper 2-join $(X_1,X_2)$ with $X=X_1$. Let $a$ and $b$ be the ends of the marker component of $T_X$.
Suppose that there is a path $P$ from $u$ to $v$ in $T$ of odd length. If $P \subseteq T|X_1$, then $P$ is a path of $T_X$, a contradiction, so $P \cap (T|X_2)$ is not empty. Since both end points of $P$ belong to
$X_1$, it follows that some edge of $P$ has one end in $X_1$ and the other
in $X_2$. By symmetry we may assume that there exist $a_1 \in A_1$
and $a_2 \in A_2$ such that $a_1a_2$ is an edge of $P$ and
$N_P(a_1) \setminus \{a_2\} \not \in A_2$.  Since $A_1$ is complete
to $A_2$, it follows that $A_2 \cap V(P)=\{a_2\}$.
Suppose that there is no  edge $b_1b_2$ of  $P$
such that $b_1 \in B_1$, $b_2 \in B_2$ and
$N_P(b_1) \setminus \{b_2\} \not \in B_2$.
Since both ends of $P$ are in $X_1$, it follows that
$V(P) \cap X_2=\{a_2\}$, and replacing $a_2$ with $a$, we obtain an
odd path from $u$ to $v$ in $T_X$, a contradiction.
This proves that such an edge $b_1b_2$ exists; now by symmetry
$B_2 \cap V(P)=\{b_2\}$.

Suppose that $N_P(a_2) \not \subseteq A_1$. Since both ends
of $P$ are in $X_1$, it follows that $P \cap (T|X_2)$ is a path
$Q$ in $T|X_2$ from $A_2$ to $B_2$.  By construction of the marker component, it follows that the path from $a$ to $b$ in the marker component of $T_X$ has the same parity as $Q$. Let $Q'$ be the path from $a$ to $b$ in the marker component of $T_X$. Now, $P' = (P \setminus Q) \cup Q'$ induces a path of $T_X$ of the same parity as $P$, a contradiction. We deduce that
$N_P(a_2)  \subseteq A_1$. Similarly, $N_P(b_2) \subseteq B_1$.
It follows that $P \setminus X_1=\{a_2,b_2\}$.
But now the path obtained from $P$ by replacing $a_2$ with $a$,
and $b_2$ with $b$ is and odd path from $u$ to $v$ in $T_X$,
a contradition.
\bbox

\subsection{Handling 2-joins and their complements} First, we show that it remains to consider trigraphs that admit proper $2$-joins.

In our next proof we will use Theorem~5.1 of \cite{Bergetrigraphs}.
The statement of that theorem require several additional definitions and
we chose not to include it here.
\begin{theorem}\label{thm:2jc}
Let $T$ be a trigraph in $\mathcal{F}$ with no balanced skew-partition and no antihole. Then, either $T$ is basic, or $T$ admits a proper 2-join.
\end{theorem}

\Proof By Theorem~\ref{lemma:t-2join}, we may assume that $\overline{T}$ admits a proper 2-join $(X_1,X_2)$ with split $(A_1, B_1,  C_1, A_2, B_2, C_2)$. By the definition of balanced skew-partition, a trigraph admits a balanced skew-partition if and only if its complement admits a balanced skew-partition. Thus, by Lemma Theorem~\ref{starcutset}, $\overline{T}$ has no star cutset.

Suppose that  $C_1 \neq \emptyset$. Let $D_1$ be an anticomponent of $C_1$. Then $D_1$
is strongly complete to $X_2 \cup (C_1 \setminus D_1)$. Let
$A_1'$ be the set of vertices in $A_1$ that are not strongly complete to
$D_1$, and define $B_1'$ similarly.
If $B_1'= \emptyset$,
then for every $a_2 \in A_2$, $A_1 \cup \{a_2\}$ is a starcutset in $\overline{T}$, a contradiction. Thus $A_1' \neq \emptyset$, and $B_1' \neq \emptyset$.

Suppose that some $a_1 \in A_1'$ has a neighbor $b_1 \in B_1'$. Then there is an
antipath $P$ of length at least two from $a_1$ to $b_1$ with interior in $D_1$.
If there is an antipath $Q$ in   $T|X_2$ with one end in $A_2$, the other end in $B_2$ and interior in $C_2$, then $T|(V(P) \cup V(Q))$  is 
an antihole of length at least five, a contradiction, so no such
$Q$ exists. Reversing the roles of $X_1$ and $X_2$, we deduce that $C_2=\emptyset$. By the definition of a 2-join and symmetry, we may assume that $|A_2|>1$.
Let $a_2 \in A_2$. Now $A_2 \cup  \{a_2\}$ is a star cutset
in $\overline{T}$ separating $A_2 \setminus a_2$ from $V(T) \setminus (A_1 \cup A_2)$, a contradiction. This proves that  $A_1'$ is strongly anticomplete to
$B_1'$.

Now $(A_1' \cup B_1', X_2 \cup (X_1 \setminus (A_1' \cup B_1')))$ is a skew-partition of $T$. We will obtain a contradiction by showing that this is
a balanced skew-partition. Suppose there is an odd path $P$ with ends
$u,v \in X_2 \cup (X_1 \setminus (A_1' \cup A_2'))$ and interior in $A_1' \cup B_1'$. We may assume that $P \setminus \{u,v\} \subseteq A_1'$. Since
$B_2 \cup C_2$ is strongly complete to $A_1'$, and $A_2$ is strongly anticomplete to $A_1'$,
it follows that $u,v \in X_1$.  Since $u \d P \d v \d a_2 \d u$ is not
an odd hole in $T$ for any $a_2 \in A_2$, it
follows that $A_1 \cap \{u,v\} \neq \emptyset$. Since $D_1$ is complete to
$A_1 \setminus A_1'$, we deduce that $v \not \in D_1$, and so
$D_1$ is strongly complete to $\{u,v\}$. Now $P$ is an odd path
with ends strongly complete to $D_1$ and with no vertex strongly complete
to $D_1$ in its interior, and every $a_2 \in A_2$ is strongly complete to
$D_1$ and strongly anticomplete to $P$, contrary to
Theorem~5.1 of \cite{Bergetrigraphs}.
Next suppose that there is an odd antipath $Q$ with ends $u,v \in A_1' \cup B_1'$ and interior in $X_2 \cup (X_1 \setminus (A_1' \cup B_1'))$.
Since
$u$ is adjacent to $v$, we may assume that $u,v \in A_1'$.
Since $u \d Q \d v \d a_2 \d u$ is not an antihole in
$T$ for any $a_2 \in A_2$, it follows that
$Q \setminus \{u,v\} \not \subseteq D_1$. Since
$Q \setminus \{u,v\}$ is anticonnected, there is an anticomponent $D_2$
of $X_2 \cup (X_1 \setminus (A_1' \cup B_1'))$ such that
$Q \setminus \{u,v\} \subseteq D_2$, and $D_2 \neq D_1$. Since $u,v \in A_1'$,
there is an antipath $Q'$ with ends  $u,v$ and interior in $D_1$.
But now $u \d Q \d v \d Q' \d v$ is an antihole in $T$, a
contradiction.

This proves that  $C_1=\emptyset$, and $C_2=\emptyset$ by symmetry. Now, $(A_1, B_1, C_1, B_2, A_2, C_2)$ is a split of a proper 2-join of $T$. This completes the proof.
\bbox

\begin{theorem}\label{thm:final}
Let $T$ be a trigraph in $\mathcal{F}$ with no balanced skew-partition, no odd prism, and no antihole. If $T$ admits a proper $2$-join, then $T$ is either complete or has an even pair disjoint from its switchable component.
\end{theorem}

\Proof Let $T$ be the vertex-minimal counterexample to the claim. Let $(A_1,B_1,C_1,A_2,B_2,C_2)$ be a split of a proper $2$-join of $T$ with $X_1=A_1\cup B_1\cup C_1$ and $X_2=A_2\cup B_2 \cup C_2$. Let $T_{X_1}$ and $T_{X_2}$ be the
blocks of the decomposition. If $T$ has a switchable component $D$, we may assume up to symmetry that $V(D)\subseteq X_2$, as no switchable pair meets both $X_1$ and $X_2$. Then, the only switchable component of $T_{X_1}$ is its marker component. By Theorem~\ref{lm:opahinfrag}, $T_{X_1}\in \mathcal{F}$, and $T_{X_1}$ admits no balanced skew-partition, no antihole, and no odd prism. 

Next, we show that $T_{X_1}$ is not complete: if $C_1$ is not empty, then a vertex $y\in C_1$ is not adjacent to any vertex of the marker component of $T_{X_1}$; if $C_1$ is empty, then by Theorem~\ref{thm:2jnobsp}, $|A_1|\geq 2$ and $T_{X_1}|A_1$ is strongly anticomplete to $b$. Now, we show that $T_{X_1}$ has an even pair disjoint from its switchable component. By Theorem~\ref{thm:2jc}, we may assume that $T_{X_1}$ is basic as $T$ is the vertex-minimal counterexample. By Theorem~\ref{thm:2jnobsp}, $|X_i|\geq 4$ for $i=1,2$, and the marker component of $T_{X_1}$ has either two or three vertices (depending on the parity of the $2$-join $(X_1,X_2)$). In either case, we have $|T_{X_1}|\geq 6$. By Theorem~\ref{switcompdom}, $T_{X_1}$ is favorable as it is not complete. Therefore, $T_{X_1}$ has an even pair $\{u,v\}$ disjoint from its marker component by Theorem~\ref{thm:basicdisj}. However, by Theorem~\ref{thm:epinfrag}, $\{u,v\}$ is also an even pair in $T$, a contradiction.
\bbox

\subsection{Proof of the main theorem} Now, we are ready to prove Theorem~\ref{thm:main}. We restate it here for the sake of clarity.

\begin{theorem}[\ref{thm:main}]
If $G$ is a Berge graph with no odd prism and no antihole of length at least six, and $G$ does not admit a balanced skew-partition, then $G$ is either complete or has an even pair.
\end{theorem}

\Proof First, $G\in \mathcal{F}$ as $G$ is Berge and has no switchable component. So the result follows from Theorems~\ref{thm:2jc}, \ref{thm:final}, and \ref{thm:basicdisj}. 
\bbox

\bibliographystyle{abbrv} 
\bibliography{refs}

\end{document}